\documentclass[12pt,a4paper]{article}
\usepackage{inputenc, amssymb, amsmath, amsthm}
\usepackage[english]{babel}
\makeatletter
\thm@style{plain}

\def\th@plain{%
  \thm@headfont{\bfseries}%
  \itshape % body font
  \thm@notefont{\rm}%
}
\def\thm@indent{\hspace*{\parindent}}
\makeatother

\def\({\left(}
\def\){\right)}

\newcommand{\ep}{\varepsilon}
\newcommand{\eps}{\varepsilon}

\newcommand{\be}{\begin{equation}}
\newcommand{\ee}{\end{equation}}

\let\epsilon\varepsilon
\let\phi\varphi
\let\le\leqslant
\let\ge\geqslant

\newtheorem{theorem}{Theorem}
\newtheorem{lemma}{Lemma}

\newcommand\inte{\int\limits}
\newcommand{\ol}[1]{ { \overline{#1} } }

\begin{document}

\centerline{\bf\uppercase{An analogue of Selberg's formula}}
\centerline{\bf\uppercase{for Motohashi's product}}

\medskip

\centerline{Sergei~N.~Preobrazhenski\u i\footnote[1]{%
{\it Preobrazhenskii Sergei Nikolayevich} ---
Department of Mathematical Analysis, Faculty of Mechanics and Mathematics,
Lomonosov Moscow State University.}}

\bigskip

\hbox to \textwidth{\hfil\parbox{0.9\textwidth}{%
\small We prove an analogue of Selberg's explicit formula for Motohashi's product
(see {\tt arXiv:1104.1358v3 [math.NT]}).
We also provide a zero-density theorem for the product, which follows from
Soundararajan's theorem for moments of the Riemann zeta-function
on the critical line.

\medskip

}\hfil}

\bigskip

%------------------------------------------------------------------

{\bf 1. Introduction.} In~\cite{Preob1} we gave a proof of the
following theorem.

Assume RH. Then there exists an infinite sequence of pairs of real numbers $(T_1,T_2)$,
$T_1=T$, $T_2=T+H$, with arbitrarily large values of $T$ and $H=c(\log\log T)^{-1}$,
such that
\[
|\zeta(1+iT_1)||\zeta(1+iT_2)|\ll(\log\log T)^{-2}
\]
and
\[
\begin{split}
&(\log\log T)^7|\zeta(1+iT_1)|^4|\zeta(1+iT_2)|^4\\
{}+&(\log\log T)^7\zeta(1+i(T_1+H))^2\zeta(1-i(T_1-H))^2\zeta(1+i(T_2+H))^2\zeta(1-i(T_2-H))^2\\
{}+&(\log\log T)^7\zeta(1+i(T_1-H))^2\zeta(1-i(T_1+H))^2\zeta(1+i(T_2-H))^2\zeta(1-i(T_2+H))^2\\
{}\ll&(\log\log T)^{-1}.
\end{split}
\]
We put
\[ %\label{PreobXzassum}
%X=\exp(0{.}5DA\log\log T\log\log\log T),\quad
z=\exp(A\log\log T\log\log\log T)
\]
with the same constant $A$ %and $D$
as in~\cite{Preob1},
set
\[
\xi(d)=\lambda_d(z)=
\begin{cases}
\mu(d)&\text{if}\quad d<z,\\
\mu(d)\frac{\log\(z^2/d\)}{\log z}&\text{if}\quad z\le d<z^2,\\
0&\text{otherwise},
\end{cases}
\]
\[
\begin{split}
&P_d(s,T_1,T_2)\\
{} &=\prod_{p\mid d}\left(1-\left(1-\frac1{p^s}\right)\left(1-\frac1{p^{s-iT_1}}\right)
\left(1-\frac1{p^{s+iT_2}}\right)\left(1-\frac1{p^{s-i(T_1-T_2)}}\right)
\left(1-\frac1{p^{2s-i(T_1-T_2)}}\right)^{-1}\right),
\end{split}
\]
and
write %$H=0$
\begin{gather*}
J(s,T_1,T_2)=\frac{\zeta(s)\zeta(s-iT_1)\zeta(s+iT_2)\zeta(s-i(T_1-T_2))}
{\zeta(2s-i(T_1-T_2))},\\
K(s,T_1,T_2)=\sum_{d\le z^2}\lambda_d(z)P_d(s,T_1,T_2).
\end{gather*}

The theorem below shows that on RH the absolute value of the product $JK$ cannot exceed
a quite small function ($(\log\log T)^{\ep}$) in a quite large neighborhood
(within the distance $c(\log\log T)^{-1}$) to the left of the line $\Re s=1$.

\begin{theorem}\label{PreobRHlargefree}
Assume the Riemann hypothesis.
Let $s_0=\sigma_0+it_0$ be a point
such that
\[ %\be\label{PreobJKassumg}
|J(s_0,T_1,T_2)K(s_0,T_1,T_2)|\ge(\log\log T)^{\ep}
\] %\ee
with arbitrarily small fixed $\ep>0$, and
\be\label{Preobsigma0t0assum}
\sigma_0=1-\frac{E_0}{\log\log T}\ge1-\frac E{\log\log T},\quad C\log\log\log T\le|t_0|\le T/2. %2H if H is larger than C\log\log\log T
\ee
Then $E_0\ge c_2(\ep)>0$.
\end{theorem}

Here we call the product $J(s_0,T_1,T_2)K(s_0,T_1,T_2)$ Motohashi's product
(see~\cite{PreobMoto}).

In this paper, we prove an analogue of Selberg's
formula for the logarithmic derivative of Motohashi's product.
We also provide a zero-density theorem for the product, which follows from
Soundararajan's theorem for moments of the Riemann zeta-function
on the critical line~\cite{PreobSound}.

{\bf 2. Lemmas and the results.}
The following lemmas are found in Ingham~\cite{PreobIngham},
theorems 26 and 27. The proofs use the functional equation
for the zeta-function.

\begin{lemma}\label{PreobInghamp71} There exists a sequence of numbers
$T_2$, $T_3$, $\ldots$, such that
\[
m<T_m<m+1\quad(m=2,3,\ldots)
\]
and
\[
\left|\frac{\zeta'}{\zeta}(s)\right|<A\log^2t\quad(-1\le\sigma\le2,t=T_m).
\]
\end{lemma}

\begin{lemma}\label{PreobInghamp73} In the region obtained by removing from the
half-plane $\sigma\le-1$ the interiors of a set of circles of radius $\frac12$ with
centres at $s=-2$, $-4$, $-6$, $\ldots$, i.e. in the region defined by
\[
\sigma\le-1,\quad|s-n|\ge\frac12\quad(n=-2,-4,-6,\ldots),
\]
we have
\[
\left|\frac{\zeta'}{\zeta}(s)\right|<A\log(|s|+1).
\]
\end{lemma}

The following lemma is a consequence of lemmas~\ref{PreobInghamp71}
and~\ref{PreobInghamp73} (cf. Lemma 9 of~\cite{PreobSelberg}).

\begin{lemma}\label{PreobSelberg9p233} There exists a sequence of numbers
$T_2$, $T_3$, $\ldots$, such that
\[
m<T_m<m+1\quad(m=2,3,\ldots)
\]
and
\[
\left|\frac{\zeta'}{\zeta}(s)\right|<A\log^2m,
\]
for $\sigma\ge-m-1/2$, $t=\pm T_m$, or $\sigma=-m-1/2$, $|t|<T_m$.
\end{lemma}

Consider the function $Z(s)=J(s,T_1,T_2)K(s,T_1,T_2)$.
By Lemma~3 of~\cite{Preob1}, for $\sigma>1$ $Z(s)$ is given
by the convergent Dirichlet series
\[
Z(s)=\sum_{n=1}^{\infty}\sigma_{iT_1}(n)\sigma_{-iT_2}(n)\left(\sum_{d\mid n}\xi(d)\right)n^{-s}.
\]
Recalling the definition of $\xi(d)$ we see that for $\sigma$
sufficiently large we can make the power series expansion for $\log Z(s)$
as $\log(1+x)$ and represent this as the convergent Dirichlet series.
Alternatively, we can expand $\log K(s)$ and use the known expansion for $\log J(s)$.
Differentiating term by term, we obtain the Dirichlet series for
$\frac{Z'}Z(s)$. Denote the coefficients of the series by $\Sigma(n)$:
\[
\frac{Z'}Z(s)=\sum_{n=2}^{\infty}\frac{\Sigma(n)}{n^s}.
\]
The following formula is an analogue of Selberg's formula for $\frac{\zeta'}{\zeta}(s)$.
We denote the complex zeros of $\zeta(s)$ by $\rho=\beta+i\gamma$,
the zeros of $K(s,T_1,T_2)$ by $r$, and the poles of $K(s,T_1,T_2)$ by $\nu$.

\begin{theorem}\label{PreobSelberg10p233} Let $x>1$ and write
\[
\Sigma_x(n)=\begin{cases}
\Sigma(n)&\text{for}\quad1\le n\le x,\\
\Sigma(n)\frac{\log^2\frac{x^3}n-2\log^2\frac{x^2}n}{2\log^2x}
&\text{for}\quad x\le n\le x^2,\\
\Sigma(n)\frac{\log^2\frac{x^3}n}{2\log^2x}
&\text{for}\quad x^2\le n\le x^3.
\end{cases}
\]
Then for
\[
\begin{split}
s\not\in&S_1=\{1,1+iT_1,1-iT_2,1+i(T_1-T_2)\},\quad s\not\in\ol{S}_1=\left\{\frac12+\frac i2(T_1-T_2)\right\},\\
s\not\in&S_{\rho}=\{\rho,\rho+iT_1,\rho-iT_2,\rho+i(T_1-T_2)\},\quad s\not\in\ol{S}_{\rho}=\left\{\frac{\rho}2+\frac i2(T_1-T_2)\right\},\\
s\not\in&S_{-2q}=\{-2q,-2q+iT_1,-2q-iT_2,-2q+i(T_1-T_2)\},\quad s\not\in\ol{S}_{-2q}=\left\{-q+\frac i2(T_1-T_2)\right\}\\
&(q=1,2,3,\ldots),\\
s\not\in&S_{\nu}=\left\{\nu=i\left(\frac{\pi k}{\log p}+\frac12(T_1-T_2)\right)\quad
\(k\in{\mathbb Z}'\subset{\mathbb Z},\quad p<z^2\)\right\},\\
s\not\in&S_r=\{r\},
\end{split}
\]
we have
\begin{gather}
\frac{Z'}Z(s)=\sum_{n<x^3}\frac{\Sigma_x(n)}{n^s}
+\sum_{u\in S_1\cup\ol{S}_{\rho}\cup\ol{S}_{-2q}\cup S_{\nu}}
\frac{x^{u-s}(1-x^{u-s})^2}{\log^2x(u-s)^3}\notag\\
\label{PreobSelbergtype}
-\sum_{u\in\ol{S}_1\cup S_{\rho}\cup S_{-2q}\cup S_r}
\frac{x^{u-s}(1-x^{u-s})^2}{\log^2x(u-s)^3},
\end{gather}
where the zeros and the poles are counted according to their multiplicities.
\end{theorem}

For the proof, we follow the argument of Selberg in Lemma 10 of~\cite{PreobSelberg}.
Note that the analogue of Lemma~\ref{PreobSelberg9p233} for the function $K(s,T_1,T_2)$
should be with the bound
\be\label{PreobdlogK}
\left|\frac{K'}K(s,T_1,T_2,z)\right|<A_{T_1,T_2,z}m,
\ee
for a sequence of contours $C_m$:
$\sigma_m\le\sigma\le\alpha$, $t=\pm T_m$, or $\sigma=\sigma_m$, $|t|<T_m$
($\sigma_m\to-\infty$, $T_m\to+\infty$). The sequence $\sigma_m$, $T_m$
is chosen properly so as to avoid zeros and poles of $K$ (and zeros of $J$).
This can be done using Jensen's theorem for $K$.
In view of the bound~\eqref{PreobdlogK} we need degree $3$ in the denominator
for the integral to tend to zero.

Next, we give a zero-density theorem for Motohashi's product.
Recall that
\[
K_X(s)=K(s,T_1,T_2)=\sum_{d\le X^2}\lambda_d(X)P_d(s,T_1,T_2),
\]
and
\[
J(s)K_X(s)=1+\sum_{n=X}^{\infty}
\sigma_{iT_1}(n)\sigma_{-iT_2}(n)
\left(\sum_{d\mid n}\lambda_d(X)\right)n^{-s}.
\]
Let
\[
f_X(s)=J(s)K_X(s)-1,
\]
and define $N_K(\sigma,T)$ to be the number of zeros $r$ of the function $K_X(s)$
such that $\Re r>\sigma\ge\frac12$, $0<\Im r\le T$.

\begin{lemma}\label{PreobTitch9dot16} If for some $X=X(\sigma,T)$, $X<T^{A_1}$,
\[
\inte_0^T|f_X(\sigma+it)|^{1/\varkappa}dt=
O\left(T^{l(\sigma)}\log^mT\right)
\]
as $T\to\infty$, uniformly for $\alpha\le\sigma\le\beta$, where $l(\sigma)$ is a positive
non-increasing function with a bounded derivative, $\varkappa$ lies between $\frac12$ and $\frac34$,
and $m$ is a constant $\ge0$, then for any fixed $A_2\ge1$
\[
N_K(\sigma,T)=O\left(T^{l(\sigma)}\log^{m+1}T\right)+O\left(\frac{T\log T}{X^{A_2}}\right)
\]
uniformly for $\alpha+1/\log T\le\sigma\le\beta$.
\end{lemma}

This lemma is a modified form of Theorem 9{.}16 of~\cite{PreobTitchmarsh}.

\begin{lemma} %\label{PreobSoundMom}
Assume RH.
For every positive real number $k$ and every $\eps>0$ we have
\[
\inte_0^T\left|\zeta\(\frac12+it\)\right|^{2k}dt\mathrel{{\ll}_{k,\eps}}
T(\log T)^{k^2+\eps}.
\]
\end{lemma}

This is a theorem of K.~Soundararajan~\cite{PreobSound}.

The following lemma is an analogue of Theorem 9{.}19 (B) of~\cite{PreobTitchmarsh}
(Ingham's zero density theorem) and follows from Soundararajan's theorem
and the upper bound for $K_X(s)$ with $X=z$ in~\cite{Preob1}.

\begin{lemma}\label{PreobTitch9dot19B} Assume RH.
%Let $T_0=\exp\(A(\log\log T)^3(\log\log\log T)^2\)$. Then
\[
\inte_0^{T}|f_X(1/2+d\log\log T\log\log\log T/\log T+it)|^{1/\varkappa}dt=
O\left(T^{1-c_d\log\log T\log\log\log T/\log T}\right).
\]
\end{lemma}

%\makeatletter
%\def\@biblabel#1{#1.}
%\makeatother

%----  ---------------------------------------

\end{document}